\documentclass[review,onefignum,onetabnum]{siamart190516}

\makeatletter
\def\myunderbrace#1{\mathop{\vtop{\m@th\ialign{##\crcr
   $\hfil\displaystyle{#1}\hfil$\crcr
   \noalign{\kern3\p@\nointerlineskip}%
   \footnotesize\downbracefill\crcr\noalign{\kern3\p@\nointerlineskip}}}}\limits}
\makeatother

\def\lesssimA#1#2{\mathrel{\vcenter{\offinterlineskip%
    \ialign{\hfil##\hfil\cr$#1<$\cr$#1\sim$\cr}%
}}}

\def\lesssim{\mathpalette\lesssimA{}}

\graphicspath{{figures/}}

\usepackage{amsmath}
\usepackage{lipsum}
\usepackage{amsfonts}
\usepackage{graphicx}
\usepackage{epstopdf}
\usepackage{todonotes}
\usepackage{algorithm}
\usepackage{algpseudocode}
\usepackage[caption=false]{subfig}
\usepackage{lineno} 

\graphicspath{{figures/}}
\ifpdf
  \DeclareGraphicsExtensions{.eps,.pdf,.png,.jpg}
\else
  \DeclareGraphicsExtensions{.eps}
\fi

\newsiamremark{remark}{Remark}
\newsiamremark{hypothesis}{Hypothesis}
\newsiamremark{example}{Example}
\crefname{hypothesis}{Hypothesis}{Hypotheses}
\newsiamthm{claim}{Claim}

\headers{Understanding the ultraspherical spectral method}{Lu Cheng and Kuan Xu}

\title{Understanding the ultraspherical spectral method}

\author{Lu Cheng\thanks{School of Mathematical Sciences, University of Science and Technology of China, 96 Jinzhai Road, Hefei, Anhui, 230026, China (\email{chengllu@mail.ustc.edu.cn},\email{kuanxu@ustc.edu.cn}).}
\and Kuan Xu\footnotemark[1]}
\usepackage{amsopn}

\nolinenumbers

\newcommand{\al}{a^{\lambda}}
\newcommand{\bg}{\boldsymbol{g}}
\newcommand{\bu}{\boldsymbol{u}}

\newcommand{\cond}{\mathrm{cond}_{E,b}}

\newcommand{\eps}{\epsilon}
\newcommand{\epsm}{\epsilon_{mach}}
\newcommand{\epsu}{\epsilon_{\min}}
\newcommand{\hR}{\hat{R}}
\newcommand{\hs}{\hat{s}}
\newcommand{\hu}{\hat{u}}
\newcommand{\lnorm}{\left \|}
\newcommand{\rnorm}{\right \|}
\newcommand{\md}{\mathrm{d}}
\newcommand{\mB}{\mathcal{B}}
\newcommand{\mD}{\mathcal{D}}

\newcommand{\mL}{\mathcal{L}}
\newcommand{\mM}{\mathcal{M}}
\newcommand{\mO}{\mathcal{O}}
\newcommand{\mP}{\mathcal{P}}
\newcommand{\mS}{\mathcal{S}}
\newcommand{\td}{\tilde}

\newcommand{\veps}{\varepsilon}

\begin{document}

\maketitle

\begin{abstract}
  The ultraspherical spectral method features high accuracy and fast solution. In this article, we determine the sources of error arising from the ultraspherical spectral method and derive its effective condition number, which explains why its backward error is consistent with a numerical method with bounded condition number. In addition, we show the cause for the Cauchy error to go below the machine epsilon and decay eventually to exact zero, revealing the fact that the Cauchy error can be misleading when used as an indicator of convergence and accuracy. The analysis in this work can be readily extended to other spectral methods, when applicable, and to the solution of PDEs.
\end{abstract}



\section{Introduction} \label{sec:intro}
The ultraspherical spectral method \cite{olv} is mathematically identical to the tau method despite the employment of polynomial spaces beyond Chebyshev. However, the ultraspherical spectral method is numerically superior to the tau method in that the former is better conditioned and can be solved in linear complexity. In this paper, we provide an error analysis for the ultraspherical spectral method by (1) identifying the sources of errors in implementing the ultraspherical spectral method, (2) giving an upper bound for the forward error in the approximate solution, (3) deriving an effective condition number, and (4) explaining the behavior of the Cauchy error and why it could be misleading.

Our discussion is structured as follows. After identifying three sources of error in \cref{sec:sources}, we revisit the Airy equation example in \cref{sec:airy}. This example motivates the analysis on the forward error in the numerical solution of the linear system that arises from the ultraspherical spectral method and the effective condition number of such a system (\cref{sec:cond}) and the behavior of the Cauchy error (\cref{sec:cauchy}). We close with a further discussion in \cref{sec:conclusion}.

\section{Three sources of error}\label{sec:sources}
Consider the ordinary differential equation
\begin{subequations}
\begin{align}
\mL u(x) &= g(x), \label{lug}\\
\mbox{s.t. } \mB  u(x) &= c, \label{buc}
\end{align}\label{intode}%
\end{subequations}
where the $N$th order linear differential operator 
\begin{align*}
\mL = a^N(x)\frac{{\md}^N}{{\md}x^N} + \ldots +a^1(x)\frac{{\md}}{{\md}x} + a^0(x) \label{opL}
\end{align*} 
for $x \in \left[-1,1\right]$ and $a^{N}(x) \neq 0$. We assume that $\mB$ contains $N$ boundary conditions and $a^0(x), \ldots, a^N(x)$ and $g(x)$ have certain regularity. In the ultraspherical spectral method, $\mL$ is expressed with the differentiation operator $\mD_{\lambda}$, the multiplication operator $\mM_{\lambda}$, and the conversion operator $\mS_{\lambda}$:
\begin{align*}
\mL = \mM_N[a^{N}]\mD_N + \sum_{\lambda = 1}^{N - 1}\mS_{N - 1}\ldots \mS_{\lambda}\mM_{\lambda}[a^{\lambda}]\mD_{\lambda} + \mS_{N - 1} \ldots \mS_{0}\mM_{0}[a^0].
\end{align*}
The simplest scenario would be $\al (x)$'s being polynomials for $0 \leq \lambda \leq N$, for which $\mM_{\lambda}[a^{\lambda}]$ and, consequently, $\mL$ are banded operators. When this is not the case, $\al (x)$ can be represented by an infinite Chebyshev series $a^{\lambda}(x) = \sum_{k = 0}^{\infty}a^{\lambda}_k T_k(x)$, resulting in $\mM_{\lambda}[a^{\lambda}]$ being dense, instead of banded. Consequently, $\mL$ is also a dense matrix of infinite dimension. Incorporating the boundary conditions by boundary bordering yields
\begin{equation}
  \begin{pmatrix}
    \mB\\
    \mL
  \end{pmatrix}\bu
  = 
  \begin{pmatrix}
    c\\
    \mS_{N-1} \ldots \mS_1\mS_0\bg
  \end{pmatrix}, \label{infode}
\end{equation}
where $\bu = (u_0, u_1, \ldots)^{\top}$ and $\bg = (g_0, g_1, \ldots)^{\top}$ are the solution vector and the Chebyshev coefficient vector of $g(x)$ respectively, both of infinite length. 

To solve \cref{infode} numerically, we have to represent everything in floating point numbers. Since it is impossible to calculate $\al_k$ accurately when $|\al_k|/\|\al(x)\|_{\infty}$ is below the machine epsilon $\epsm$, the trailing coefficients of $\mO(\epsm)$ are discarded following, for instance, the chopping strategy \cite{aur}. Thus, $\al(x)$ is replaced by a finite Chebyshev series, denoted by $\hat{a}^{\lambda}(x)$, whose coefficients are all floating point numbers. We then construct the multiplication operator using $\hat{a}^{\lambda}(x)$ and denote it by $\hat{\mM}_{\lambda}$. Similarly, we let $\hat{\mD}_{\lambda}$ and $\hat{\mS}_{\lambda}$ be the floating point representations of $\mD_{\lambda}$ and $\mS_{\lambda}$ respectively to have
\begin{align*}
\hat{\mL} = fl\left( \hat{\mM}_N[\hat{a}^N]\hat{\mD}_N + \sum_{\lambda = 1}^{N - 1}\hat{\mS}_{N - 1}\ldots \hat{\mS}_{\lambda}\hat{\mM}_{\lambda}[\hat{a}^{\lambda}]\hat{\mD}_{\lambda} + \hat{\mS}_{N - 1} \ldots \hat{\mS}_{0}\hat{\mM}_{0}[\hat{a}^0]\right),
\end{align*}
where $fl(\cdot)$ is the function that evaluate an expression in floating point arithmetic. Now $\hat{\mL}$ is a banded matrix of infinite dimension. Replacing other parts of \cref{infode} by their floating point representations likewise yields
\begin{align}
\begin{pmatrix}
    \hat{\mB}\\
    \hat{\mL}
  \end{pmatrix}\hat{\bu}
  = 
  \begin{pmatrix}
    \hat{c}\\
    \hat{\mathbf{r}}
  \end{pmatrix}, \label{infodefl}
\end{align}
where $\hat{\mathbf{r}} = fl(\hat{\mS}_{N-1} \ldots \hat{\mS}_1\hat{\mS}_0 \hat{\bg})$. Here, $\hat{\bg}$ is obtained by replacing the discarded elements in $\bg$ by zeros and those retained by their floating point representations. Hence, $\hat{\mathbf{r}}$ has only a finite number of nonzero elements. Note that $\hat{\bu} \neq fl(\bu)$ in general. Though the entrywise error introduced to $\mL$ and the right-hand side by the floating point arithmetic is in the order of $\epsm$ in a relative sense, it could lead to a disastrous error in the numerical solution $\hu$ when the problem is ill-conditioned or ill-posed. See \cite[\S 3]{tre2} for a brief discussion. If we denote this error by $\veps_F$ and let
\begin{align*}
\Phi_n = \left(T_0(x)|T_1(x)| \ldots |T_{n-1}(x) \right)
\end{align*}
be the quasimatrix formed by the first $n$ Chebyshev polynomials,
\begin{align}
  \veps_F = u(x) - \Phi_{\infty} \hat{\bu}. \label{ef}
\end{align}

Next, we make \cref{infodefl} finite dimensional with the truncation operator $\mP_n = (I_n , \boldsymbol{0})$ to obtain the $n \times n$ almost-banded system
\begin{align}
A u = f, \label{finode}
\end{align}
where
\begin{align*}
A = \mP_n 
\begin{pmatrix}
  \hat{\mB}\\
  \hat{\mL}
\end{pmatrix}\mP_n^{\top},~~
u = \mP_n \hat{\bu},~~
f = 
\begin{pmatrix}
\hat{c}\\
\mP_{n-N}\hat{\mathbf{r}}
\end{pmatrix}.
\end{align*}
It is \cref{finode} that we solve for the $n$-vector $u$, which is a finite length approximation to $\hat{\bu}$. Here note that $u$ and $u(x)$ are different, as the latter is the exact solution. Thus, the error in $u$ caused by the truncation of the system is
\begin{align}
\veps_T = \Phi_{\infty} \hat{\bu} - \Phi_n u. \label{et}
\end{align}
It is easy to see that $\veps_T$ diminishes as the degrees of freedom $n$ increases.

When solving \cref{finode} in floating point arithmetic, the best we can hope is to obtain an approximate solution $\hu$ with an error $\veps_S = u - \hu$ introduced in the course of solution due to rounding. The approximate solution can be deemed as the exact solution to a perturbed problem
\begin{align}
\left(A + \Delta A\right) \hu = f + \Delta f, \label{advf}
\end{align}
where $|\Delta A| \leq \eps E$ and $|\Delta f| \leq \eps b$ for some $\eps > 0$, $E \geq 0$, and $b \geq 0$. In addition, the error introduced to $\Phi_n u$ is
\begin{align}
\veps_R = \Phi_n u - \Phi_n \hu = \Phi_n \veps_S. \label{es}
\end{align}

So far, we have encountered four linear systems. The first, \cref{infode}, is an equivalent representation of the original problem \cref{intode} in the polynomial space spanned by $C^{(N)}$. The second, \cref{infodefl}, is the finite precision approximation of \cref{infode}. The third, \cref{finode}, is the finite-dimensional approximation of \cref{infodefl} that we aim to solve for an approximate solution. Lastly, \cref{advf} is the linear system that the approximate solution we obtain actually solves. The total error $\veps = \veps_F+\veps_T+\veps_R = u(x)-\Phi_n \hu$ by \cref{ef,et,es}. Our focus in the next two sections is $\veps_R$ or, more precisely, $\veps_S$ and its relation to the backward errors $\Delta A$ and $\Delta f$. We will also discuss the role $\veps_F$ plays in $\veps$ at the end of \cref{sec:cauchy}. 



\section{The Airy equation revisited}\label{sec:airy}
We proceed our discussion by a revisit to the Airy equation \cite[\S 3.3]{olv}
\begin{align}
\mu u'' - xu = 0, ~\mbox{s.t. }  u(-1) = \operatorname{Ai}\left(-\sqrt[3]{\frac{1}{\mu}}\right), ~ u(1) = \operatorname{Ai}\left(\sqrt[3]{\frac{1}{\mu}}\right), \label{airy}
\end{align}
where $\operatorname{Ai}(x)$ is the Airy function of the first kind, and the exact solution to \cref{airy} is the Airy function with a rescaled argument, i.e.,
\begin{align}
u(x) = \operatorname{Ai}\left(\sqrt[3]{\frac{1}{\mu}} x \right). \label{asolu}
\end{align}

\begin{figure}[tbhp]
\centering
\subfloat[condition number]{\label{fig:cond}\includegraphics[scale=0.46]{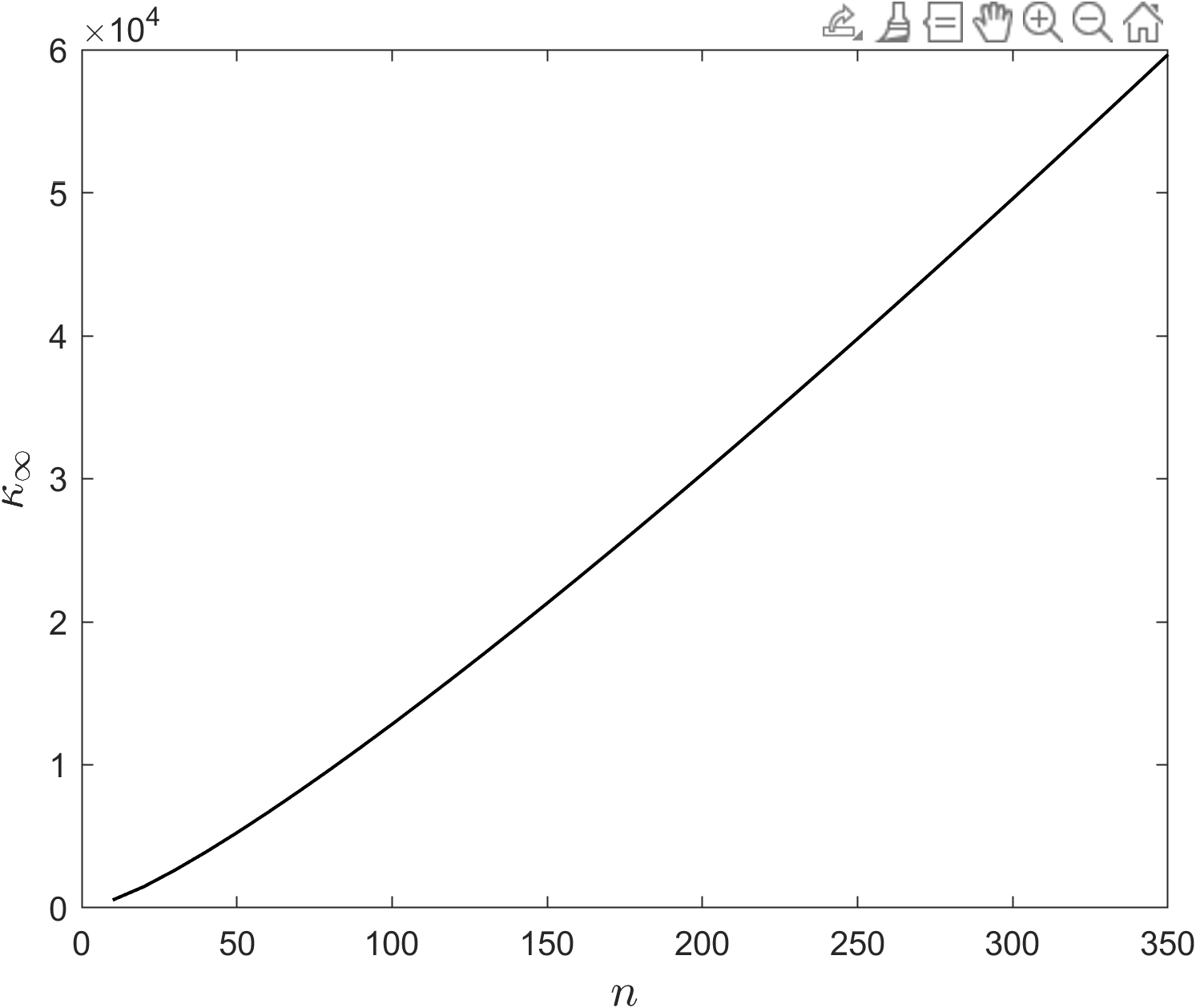}}
\hfill
\subfloat[total error]{\label{fig:f64}\includegraphics[scale=0.46]{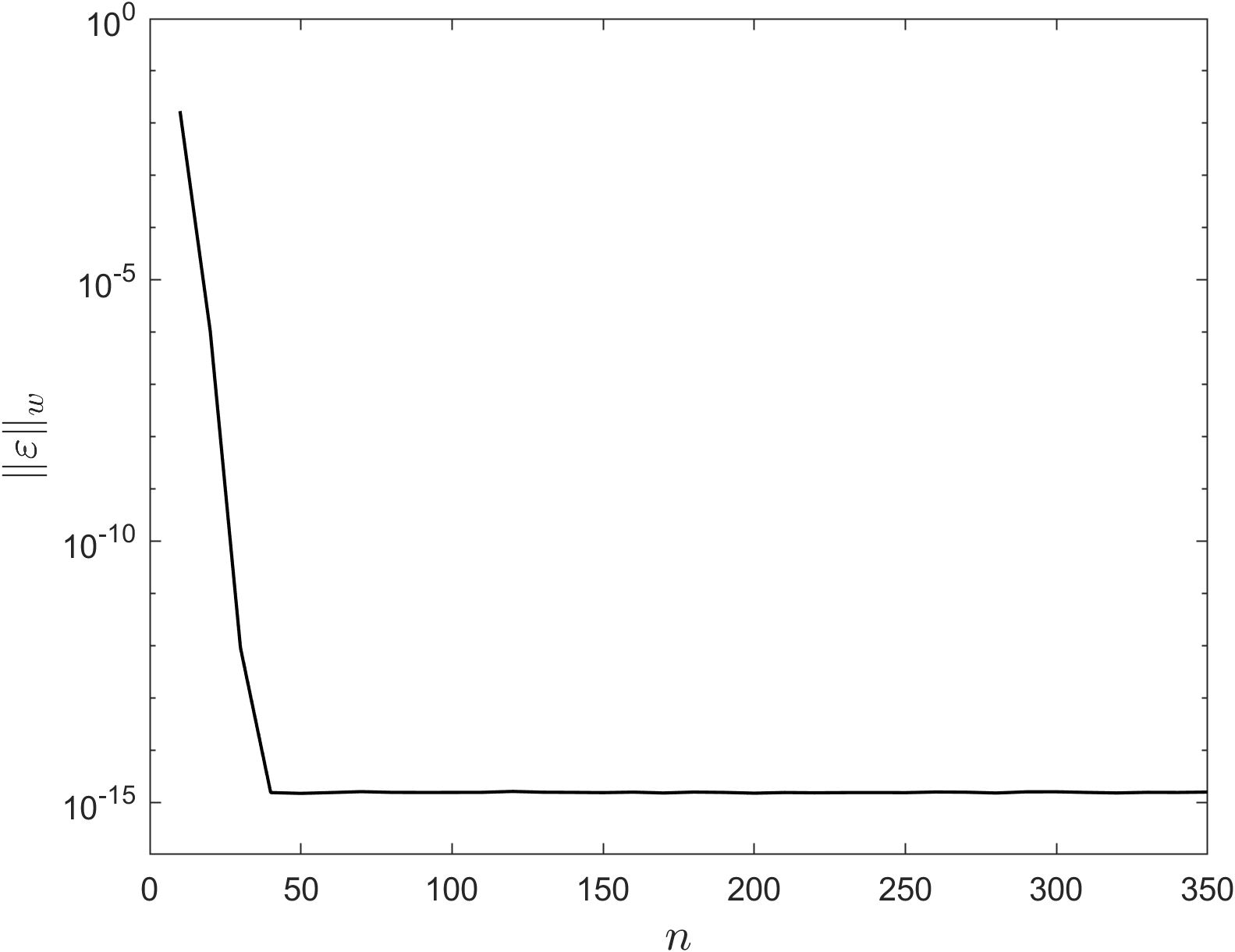}}\\
\subfloat[rounding error]{\label{fig:f64bf}\includegraphics[scale=0.46]{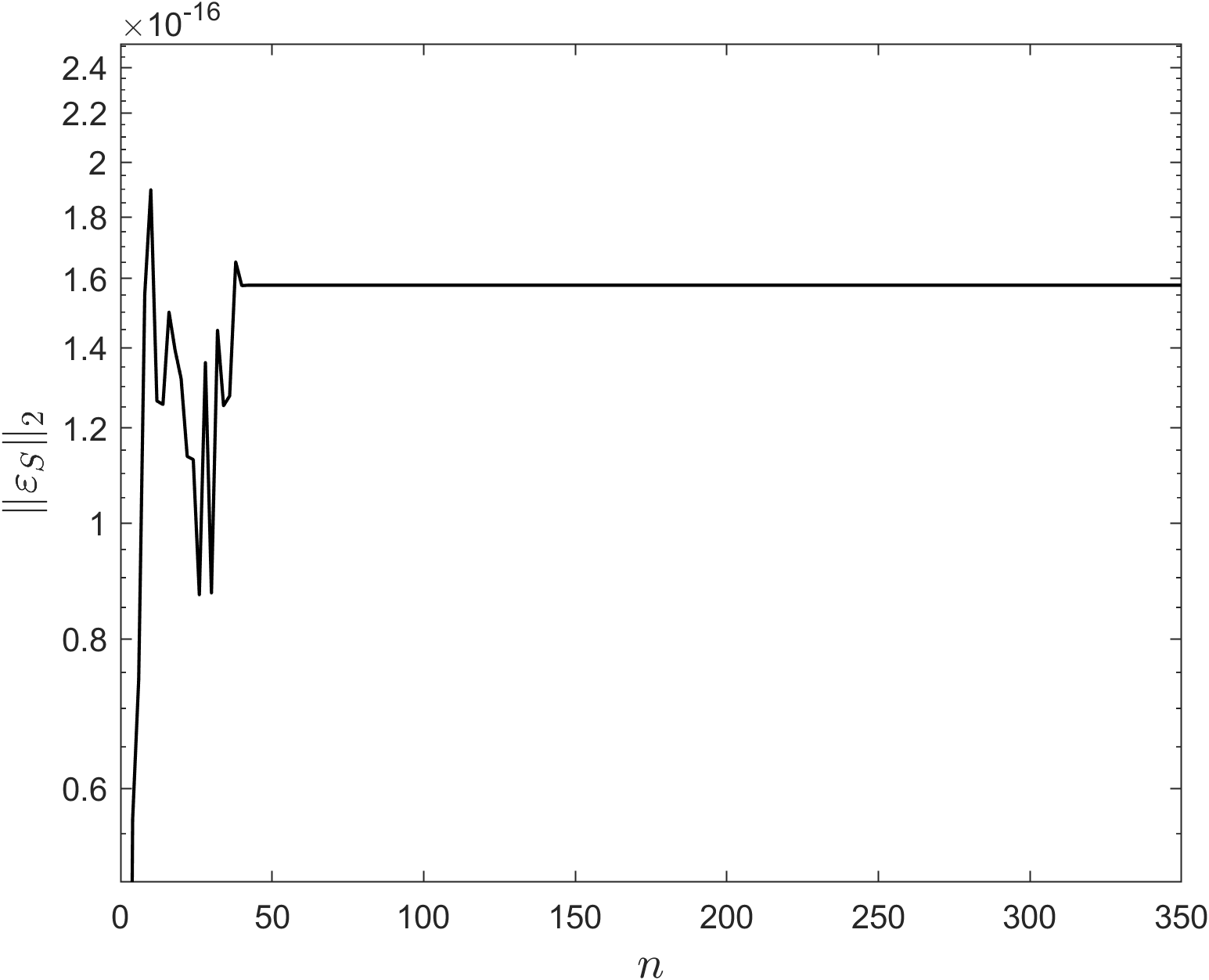}}
\hfill
\subfloat[Cauchy error]{\label{fig:cauchy}\includegraphics[scale=0.46]{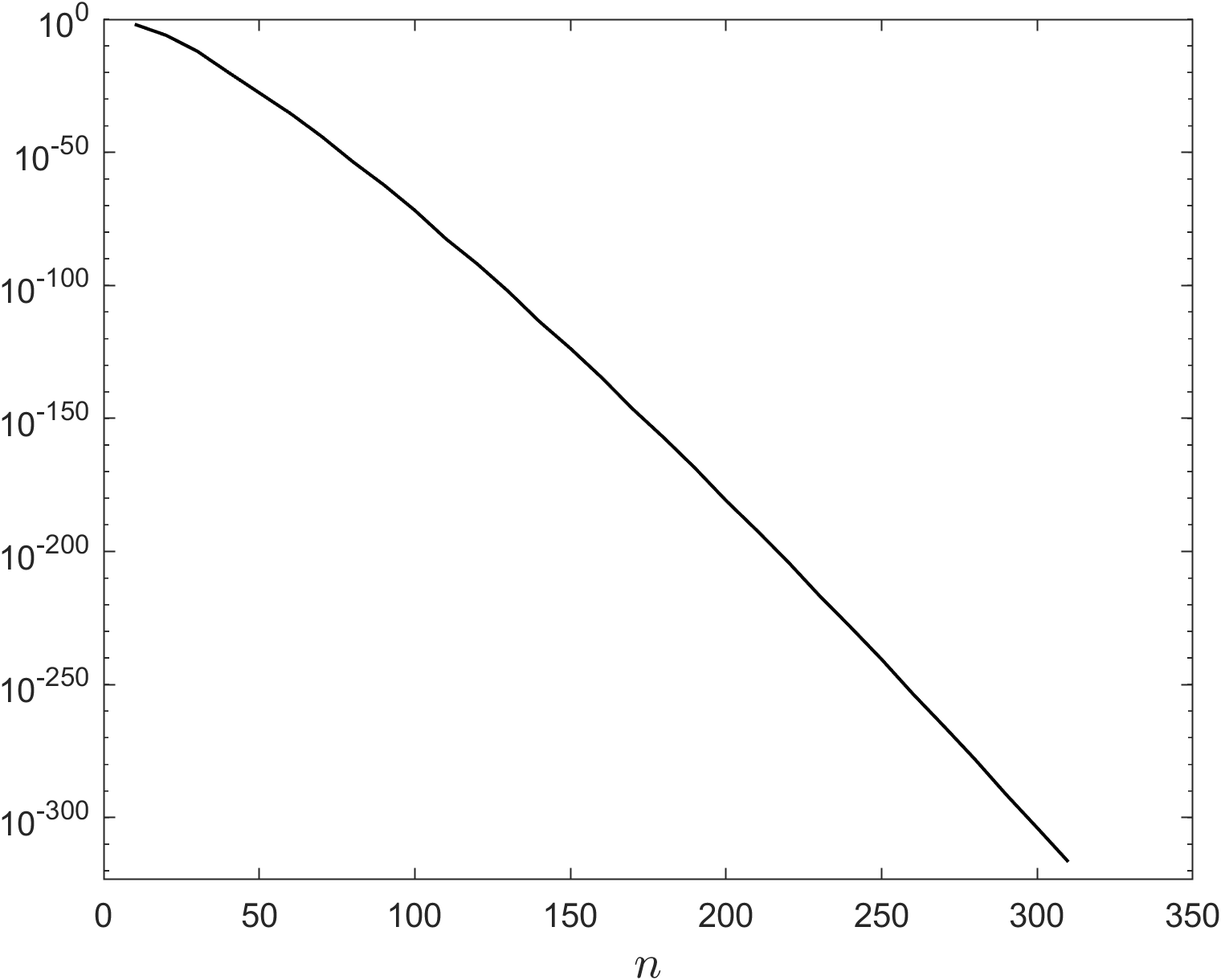}}
\caption{Solving \cref{airy} for $\mu=10^{-2}$ using the ultraspherical spectral method.}\label{fig:airy1}
\end{figure}

We solve \cref{airy} for $\mu = 10^{-2}$ with the ultraspherical spectral method in double precision and show the results in \cref{fig:airy1}. In \cref{fig:cond}, we display the condition number $\kappa_{\infty}(A) = \|A\|_{\infty}\|A^{-1}\|_{\infty}$, which grows at least linearly. 
\cref{fig:f64} is the plot of $\|\veps\|_{w}$ with $w = 1 / \sqrt{1 - x^2}$, i.e., the $2$-norm of the Chebyshev coefficients of $\veps$. It is shown that $\|\veps\|_{w}$ decays exponentially to the machine epsilon until it levels off. Here, we evaluate \cref{asolu} with \texttt{airyai} from \textsc{Julia}'s \texttt{SpecialFunctions} package in octuple precision\footnote{The $256$-bit octuple precision is the default format of \textsc{Julia}'s \texttt{BigFloat} type of floating point number. \texttt{BigFloat}, based on the GNU MPFR library, is the arbitrary precision floating point number type in \textsc{Julia}.} and use the result as the exact solution. To calculate the error, the solution in double precision is promoted to octuple precision before the difference from the ``exact'' is taken. In case of the collocation-based pseudospectral method, the error curve would bear the same fast decay but followed by a gradual yet everlasting rebound, which matches the growth of the condition number with $n$. However, the plateau in \cref{fig:f64} contradicts the growth of the condition number shown in \cref{fig:cond}.

To separate $\veps_R$ from $\veps_F$ and $\veps_T$, we also solve \cref{airy} by the ultraspherical spectral method but in octuple precision and calculate the difference between the double- and octuple-precision solutions by promoting the former to the octuple precision. This difference can be deemed as $\veps_S$ and is shown in \cref{fig:f64bf}. The plateau beyond $40$ further suggests a bounded condition number. 

\cref{fig:f64,fig:f64bf} echo the description in \cite{olv}, here we quote---\textit{the backward error is consistent with a numerical method with bounded condition number.} To explain this, Olver and Townsend take a detour---they show that \cref{finode} can be right preconditioned with a diagonal preconditioner\footnote{Note that this right preconditioner is used only for transforming the original problem to a new one whose condition number is bounded. It does not mean that applying this preconditioner can improve the accuracy whatsoever. Any additional accuracy gained in solving the preconditioned system $AD \nu= f$ would be lost altogether in recovering $u$ from $D^{-1}u = \nu$.} $D$ so that the condition number of $AD$ is bounded for any $n$, that is, $\kappa(AD) = \mO(1)$. Since solving the unpreconditioned system \cref{finode} with QR factorization enjoys the same stability and, consequently, the same condition number as the preconditioned one, they arrive at the assertion above. However, the condition number is supposed to be a method-independent quantity. We therefore wonder if a direct route can be taken to show that the condition number of $A$ is bounded without restricting ourselves to a specific method. This motivates the analysis given in \cref{sec:cond}.

In \cite{olv}, Cauchy errors are used to measure the convergence and the accuracy of the computed solution. In a few numerical experiments\footnote{See Figures 2.2 and 5.3 in \cite{olv}} of \cite{olv}, the Cauchy error keeps decaying even beyond the machine epsilon until it is out of the range of floating point numbers in double precision. Such a behavior is reproduced in \cref{fig:cauchy} for solving \cref{airy}. Specifically, we solve \cref{finode} and an $n^{\dagger} \times n^{\dagger}$ system $A^{\dagger} u^{\dagger} = f^{\dagger}$ for $n^{\dagger} = \lceil 1.01 n \rceil$,
\begin{align*}
A^{\dagger} = \mP_{n^{\dagger}} 
\begin{pmatrix}
  \hat{\mB}\\
  \hat{\mL}
\end{pmatrix}\mP_{n^{\dagger}}^{\top},~~
u^{\dagger} = \mP_{n^{\dagger}} \hat{\bu},~~
f^{\dagger} = 
\begin{pmatrix}
  \hat{c}\\
\mP_{n^{\dagger}-N}\hat{\mS}_{N-1}\hat{\mS}_{N-2} \ldots \hat{\mS}_1\hat{\mS}_0\hat{\bg}
\end{pmatrix}.
\end{align*}
The Cauchy error is then obtained as $\|\hu - \hu^{\dagger}(1:n)\|_2$, where $\hu^{\dagger}$ is the computed solution to the $n^{\dagger} \times n^{\dagger}$ system. It can be seen that Cauchy error keeps decaying as $n$ increases until it underflows to exact zero. Since the computation is done in double precision, we are supposed to be unable to calculate any quantity smaller than $\mO(10^{-16})$ faithfully using $\mO(1)$ input data. Thus, a question immediately suggests itself---why does the Cauchy error decays even beyond $\mO(10^{-16})$? This question is answered in \cref{sec:cauchy}.


\section{Conditioning of the ultraspherical spectral method}\label{sec:cond}
In this section, we show that the conditioning of \cref{finode} depends solely on certain effective parts of $A$, $A^{-1}$, and $f$. To facilitate the discussion, we use subscripts in this section to indicate the partitions of matrices and vectors. The $n$-vectors $u$, $\hu$, and $f$ are partitioned as
\begin{align*}
u = 
\begin{pmatrix}
u_1\\
u_2
\end{pmatrix},~
\hu = 
\begin{pmatrix}
\hu_1\\
\hu_2
\end{pmatrix},\text{ and }
f = 
\begin{pmatrix}
f_1\\
f_2
\end{pmatrix},
\end{align*}
where $u_1$, $\hu_1$, and $f_1$ are $k$-vectors. The length $k$ is chosen so that
\begin{subequations}
\begin{align}
f_2 = 0, \label{f2} \\
\hu_2 = 0. \label{u2}
\end{align}
Since only the first finite number of elements of $\hat{\mathbf{r}}$ are nonzero, \cref{f2} must hold when $n$ is large enough. Similarly, \cref{u2} is guaranteed by the assumption that the solution has certain regularity \cite[\S 7 \& \S 8]{tre} and the employment of the chopping algorithm for the solution $\hu$. Because of \cref{f2}, we also partition $b = (b_1^T | b_2^T)^T$ conformingly with $b_1 \in \mathbb{R}^k$ and assume
\begin{align}
b_2 = 0. \label{b2}
\end{align}\label{ufb2}
\end{subequations}

\begin{figure}[tbhp]
\centering
\subfloat[]{\label{fig:E}\includegraphics[scale=0.32]{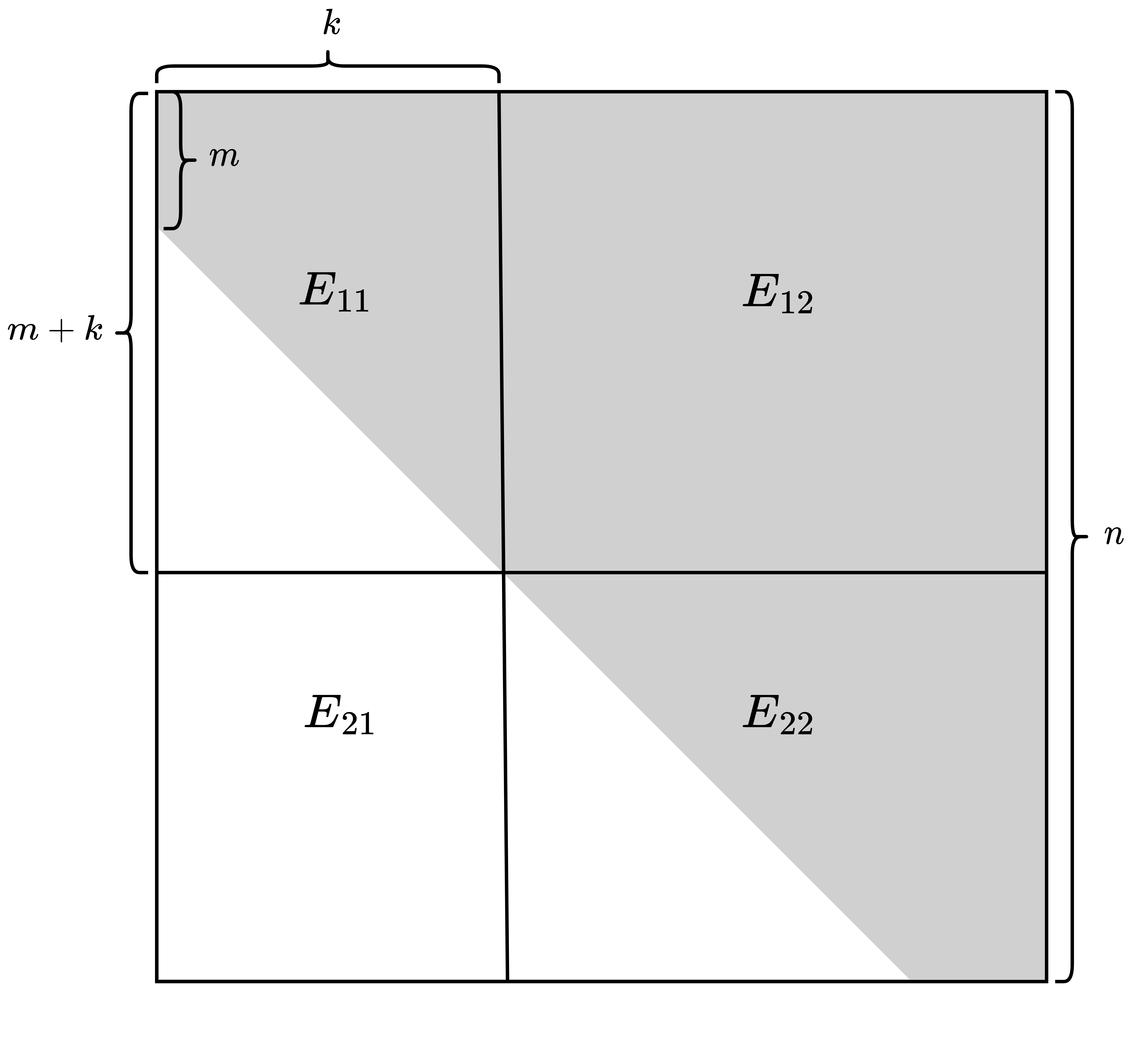}}
\hfill
\subfloat[]{\label{fig:errbound}\includegraphics[scale=0.46]{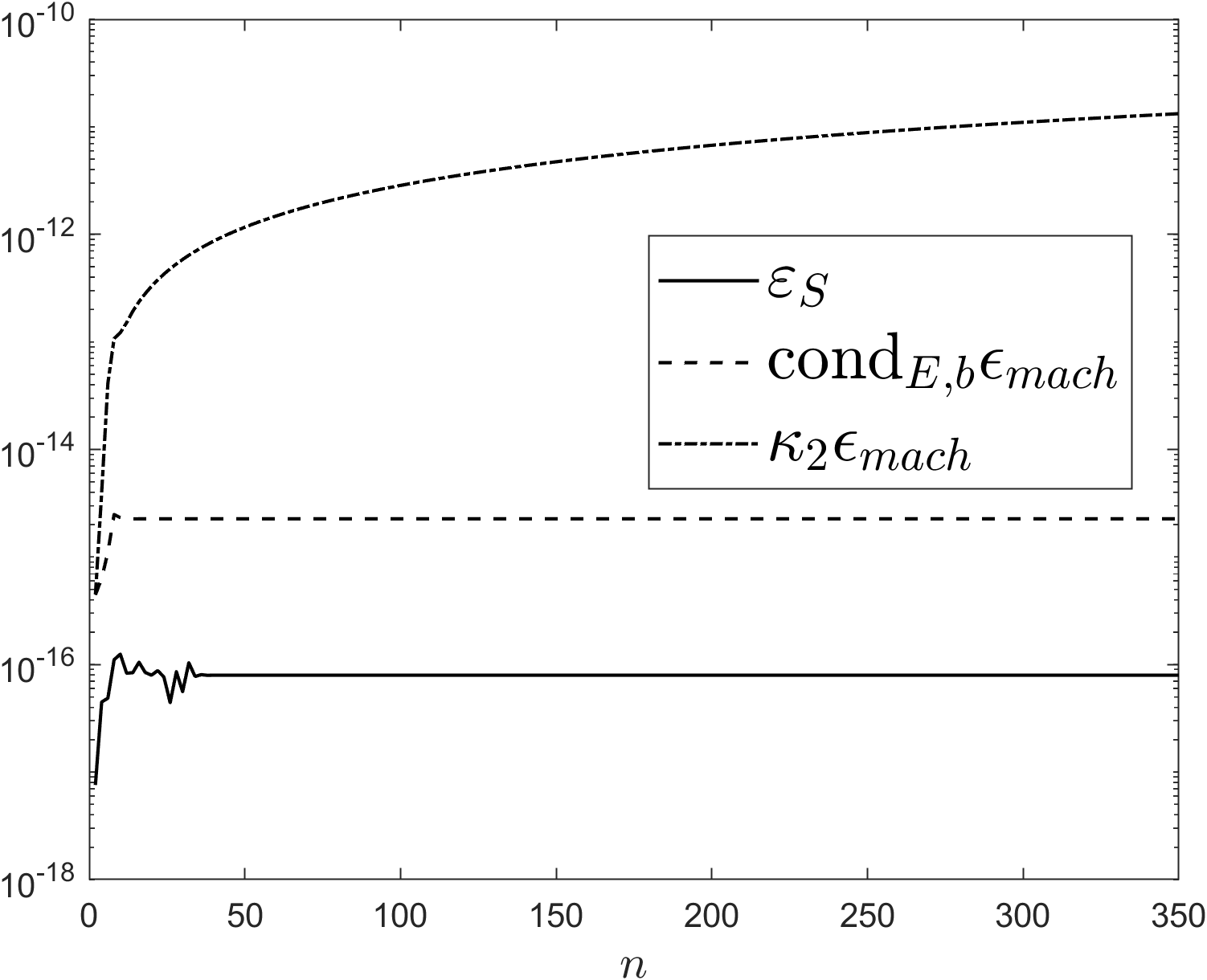}}
\caption{(a) Partition of $E$. The shade highlights the region where the entries are nonzero. (b) Estimates of the forward error in the computed solution to \cref{airy} for $\mu=10^{-2}$, compared with the actual rounding error $\veps_S$.}\label{fig:sec2}
\end{figure}

The following proposition bounds the forward error $\veps_S$ in terms of the backward error. When \cref{finode} is solved by a direct method, e.g., QR or LU, the zeros below $m$th subdiagonal remain intact. Thus, it is natural to let $E$ be $m$-Hessenberg \cite[\S 4.5]{wat}, even though $A$ is sparser as almost-banded.

\begin{proposition} \label{thm:vmu}
If $\eps \lnorm |A^{-1}|E\rnorm \leq c < 1$, where $c$ is a constant, $E$ is an $n\times n$ $m$-Hessenberg matrix, and $\lnorm \cdot \rnorm$ is an absolute norm, it follows from \cref{ufb2} that                                          
\begin{align}
\begin{split}
\|u - \hu \| \leq \frac{\epsilon}{1 - \epsilon \| |{B}| E\|} 
&\left\| |B_1| b_1 + |\tilde{B}_1|E_{11}|u_1| \right\|,
\end{split} \label{errc}
\end{align}
where $B = A^{-1} = \begin{pmatrix} {B_{1}}& {B_{2}} \end{pmatrix} = \begin{pmatrix} \td{B}_{1}& \td{B}_{2} \end{pmatrix}$ and 
$E = \begin{pmatrix} 
{E_{11}}& {E_{12}}\\
{E_{21}}& {E_{22}}
\end{pmatrix}
$ with $B_{1}\in \mathbb{R}^{n\times k}$, $\td{B}_{1} \in \mathbb{R}^{n\times (k + m)}$, and $E_{11} \in \mathbb{R}^{(k + m)\times k}$. 
\end{proposition}
\begin{proof}
Simply taking the difference of \cref{advf} and \cref{finode} to yield
\begin{align}
|u - \hu| = |B\Delta f - B\Delta A u + B\Delta A(u - \hu)| \leq \eps |B|b + \eps |B|E|u| + \eps|B|E|u - \hu|. \label{u-v}
\end{align}
By \cref{b2} and $B$'s partition, we have
\begin{align}
|B|b = |B_1| b_1. \label{Bb}
\end{align}
The partition of $E$ and the assumption that $E$ is an $m$-Hessenberg matrix imply $E_{21} = 0$, as shown in \cref{fig:E}. Thus, 
\begin{align}
|B|E|u| = |\td{B}_{1}| E_{11}|u_1|, \label{BEu}
\end{align}
where \cref{u2} is used. The inequality \cref{errc} is obtained by combining \cref{u-v,Bb,BEu}.
\end{proof}

The upper bound \cref{errc} for the forward error does not depend on $\tilde{B}_2$, i.e., the last $n-m-k$ columns of $A^{-1}$. This differs from the result for a general linear system, e.g., Theorem 7.4 in \cite[\S 7.2]{hig}, where the entire $A^{-1}$ has impact on the bound for the forward error. In addition, only the top left submatrix $E_{11}$ and the first $k$ elements of $f$ figure in the bound, whereas the traditional componentwise analysis involves the entire $E$ and $f$.

Following the definition of the componentwise condition number \cite[Eqn. $(7.11)$]{hig}, also known as Skeel's condition number, we obtain the componentwise condition number in the current context from \cref{u-v} and \cref{Bb} as 
\begin{align}
\cond = \frac{\lnorm|B_{1}| b_1 + |\td{B}_{1}| E_{11}|u_1|\rnorm_{\infty}}{\|u\|_{\infty}}.\label{cond}
\end{align}


Again, $\cond$ involves $E_{11}$, $f_1$, and the first $k+m$ columns of $A^{-1}$ only. The exclusion of other parts of $E_{11}$, $f_1$, and $A^{-1}$ suggests that $\cond$ becomes a constant for a problem so long as $n\geq k$. With $\cond$, we can estimate the forward error using the rule of thumb \cite[\S 1.6]{hig}
\begin{align}
\| \veps_S \| = \|u - \hu\| \lesssim \cond \epsm. \label{uvce}
\end{align}
This estimate is shown in \cref{fig:errbound}, where $E=|A|$ and $b=|f|$ are used for evaluating $\cond$ via \cref{cond}. For comparison, we also include the actual $\veps_S$ plotted in \cref{fig:f64bf} and the estimate obtained using \cref{uvce} but with $\cond$ replaced by $\kappa_2 = \lnorm A \rnorm_2 \lnorm A^{-1} \rnorm_2$. Apparently, $\cond$ gives a much better estimate than $\kappa_2$---not only being constant after $n=k$ but also a few orders closer to the actual $\veps_S$.

\section{The Cauchy error}\label{sec:cauchy}
Now we show why the Cauchy error decays beyond the machine epsilon until it finally becomes exactly zero due to underflow. Our analysis does not distinguish the Givens QR and the Householder QR. In this section, subscripts are used to indicate the dimension of vectors and matrices. For example, $f_n$ and $A_n$ are an $n$-vector and an $n\times n$ matrix respectively. Henceforth, we rewrite \cref{finode} as
\begin{align}
A_n u_n = f_n. \label{Auf}
\end{align}
As above, we assume again that only the first $k$ entries in $f_n$ are nonzero.

Solving \cref{Auf} with QR involves three steps---calculating the QR factorization of $A_n = Q_nR_n$, application of $Q_n^T$ to $f_n$ to form the transformed right-hand side $s_n = Q_n^Tf_n$, and solving the upper triangular system $R_nu_n=s_n$ by back substitution. Let $\hat{R}_n$ and $Q_n$ be the computed upper triangular and the exact orthonormal $QR$ factors of $A_n$ respectively such that
\begin{align*}
A_n+\overline{\Delta A}_n = Q_n \hR_n,
\end{align*}
where $\overline{\Delta A}_n$ is the backward error caused by the QR factorization only and, therefore, is different from $\Delta A$. The lemma that follows demonstrates that under a reasonable and mild assumption the sum of the squares of the last $m$ entries in any row of $Q_n$ decays exponentially.
\begin{lemma}\label{lem:q}
Suppose that $n \geq m$ and
\begin{align}
\frac{\|Q_n(j,n-m+1:n)\|_2}{\|Q_n(j,n-m:n)\|_2} \leq \sigma < 1 \label{Qc} 
\end{align}
for $j \leq n$. Then
\begin{align}
\lim_{n \to \infty}\|Q_n(j,n-m+1:n)\|_2 &\leq \lim_{n \to \infty} \sigma^{n-m}\|Q_m(j,1:m)\|_2 = 0. \label{limq}
\end{align}
\end{lemma}
\begin{proof}
Since $A_n$ and $A_{n+1}$ are $m$-Hessenberg, their respective QR factorizations differ only in the $(n-m+1)$th row and onward. That is,
\begin{subequations}
\begin{align} 
\hR_n(1:n-m,j) &= \hR_{n+1}(1:n-m,j), \label{hRn}\\
Q_n(j,1:n-m) &= Q_{n+1}(j,1:n-m). \label{QQ}
\end{align}
\end{subequations}
Since $Q_n$ and $Q_{n+1}$ are orthogonal,
\begin{align}
\|Q_n(j,:)\|_2 = \|Q_{n+1}(j,:)\|_2 = 1,
\end{align}
which, along with \cref{QQ}, gives
\begin{align*}
\|Q_n(j,n-m+1:n)\|_2 = \|Q_{n+1}(j,n-m+1:n+1)\|_2.
\end{align*}
Thus, the assumption \cref{Qc} leads to
\begin{align*}
\|Q_n(j,n-m+1:n)\|_2 \leq \sigma^{n-m} \|Q_{m}(j,1:m)\|_2,
\end{align*}
which implies \cref{limq}.
\end{proof}

Let the computed transformed right-hand side be $\hs_n = Q_n^T(f_n+\Delta f_n)$, where the backward error
\begin{align}
\|\Delta f_n\|_2 \leq \gamma \| f_n \|_2 \label{dq}
\end{align}
and $\gamma = \tilde{\gamma}_{n^2}$ for the Householder QR and $\gamma = \tilde{\gamma}_{2n-2}$ for the Givens QR \cite[\S 19.3 \& \S 19.6]{hig}. Here $\tilde{\gamma}_{n} = cn\epsm / (1 - cn\epsm)$ for a small integer constant $c$ \cite[\S 3.4]{hig}. Consider the partially nested systems 
\begin{subequations}
\begin{align}
\hR_n v_n &= \hs_n, \label{rvs1} \\
\hR_{n^+} v_{n^+} &= \hs_{n^+}, \label{rvs2}
\end{align}\label{rvs}%
\end{subequations}
where $n^+>n$. The following theorem shows that when $n$ is large enough, the computed solution $\hu_n$ and $\hu_{n^+}(1:n)$ to \cref{rvs} are identical and the rest of $\hu_{n^+}$ are zero.

\begin{theorem} \label{lem:u2e0}
Suppose that the assumption in \cref{lem:q} holds. When $n$ is sufficiently large,
\begin{subequations}
\begin{align}
\hu_n(1:n-m) &= \hu_{n^+}(1:n-m), \label{v1} \\
\hu_n(n-m+1:n) &= 0, ~\hu_{n^+}(n-m+1:n^+) = 0, \label{v2}
\end{align} \label{v}%
\end{subequations}
for any $n^+ > n$.
\end{theorem}
\begin{proof}
Equations \cref{limq} and \cref{dq} and the fact that only the first $k$ entries of $f$ are nonzero guarantee that $\|\hs_n(n-m+1:n)\|_2 < \epsu$, where $\epsu \approx 4.94 \times 10^{-324}$ is the smallest positive subnormal floating point number in double precision. Because of underflow, $\hs_n(n-m+1:n) = 0$. Since $\hR_n$ is upper triangular, 
\begin{subequations}
\begin{align}
\hu_n(n-m+1:n) = 0.
\end{align}
Because of $n^+ > n$, starting from at least the $(n-m+1)$th entry all subsequent entries of $\hu_{n^+}$ are zero, i.e., 
\begin{align}
\hu_{n^+}(n-m+1:{n^+}) = 0.
\end{align} \label{v0}%
\end{subequations}
Due to \cref{v0}, \cref{rvs} reduces to
\begin{align*}
\hR_n(1:n-m,1:n-m)v_n(1:n-m) &= \hs_n(1:n-m),\\
\hR_{n^+}(1:n-m,1:n-m)v_{n^+}(1:n-m) &= \hs_{n^+}(1:n-m).
\end{align*}
The fact that \cref{QQ} implies that
\begin{align*}
\hs_n(1:n-m) = \hs_{n^+}(1:n-m),
\end{align*}
which, along with \cref{hRn}, suggests \cref{v1}.
\end{proof}

\cref{lem:u2e0} demonstrates that the computed solutions to \cref{Auf} with different $n$ converges as their dimensions increase. Beyond a certain point, the difference between the solutions for different $n$ becomes exactly zero and remains so. Thus, we have answered the question posed at the end of \cref{sec:airy}. Apparently, the decay beyond $\epsm$ is the consequence of repeated augmentation of the system and the $Q$ factor whose entries have decaying magnitudes as its dimension grows, therefore being merely an \emph{artifact}. Hence, the decay of Cauchy errors beyond $\epsm$ should be ignored completely.

\begin{figure}[tbhp]
  \centering
  \subfloat[Cauchy error and $\sum_{j=1}^k\|Q_n(j,n-m+1:n)\|_2$ versus $n$.]{\label{fig:CauchyQ}\includegraphics[scale=0.46]{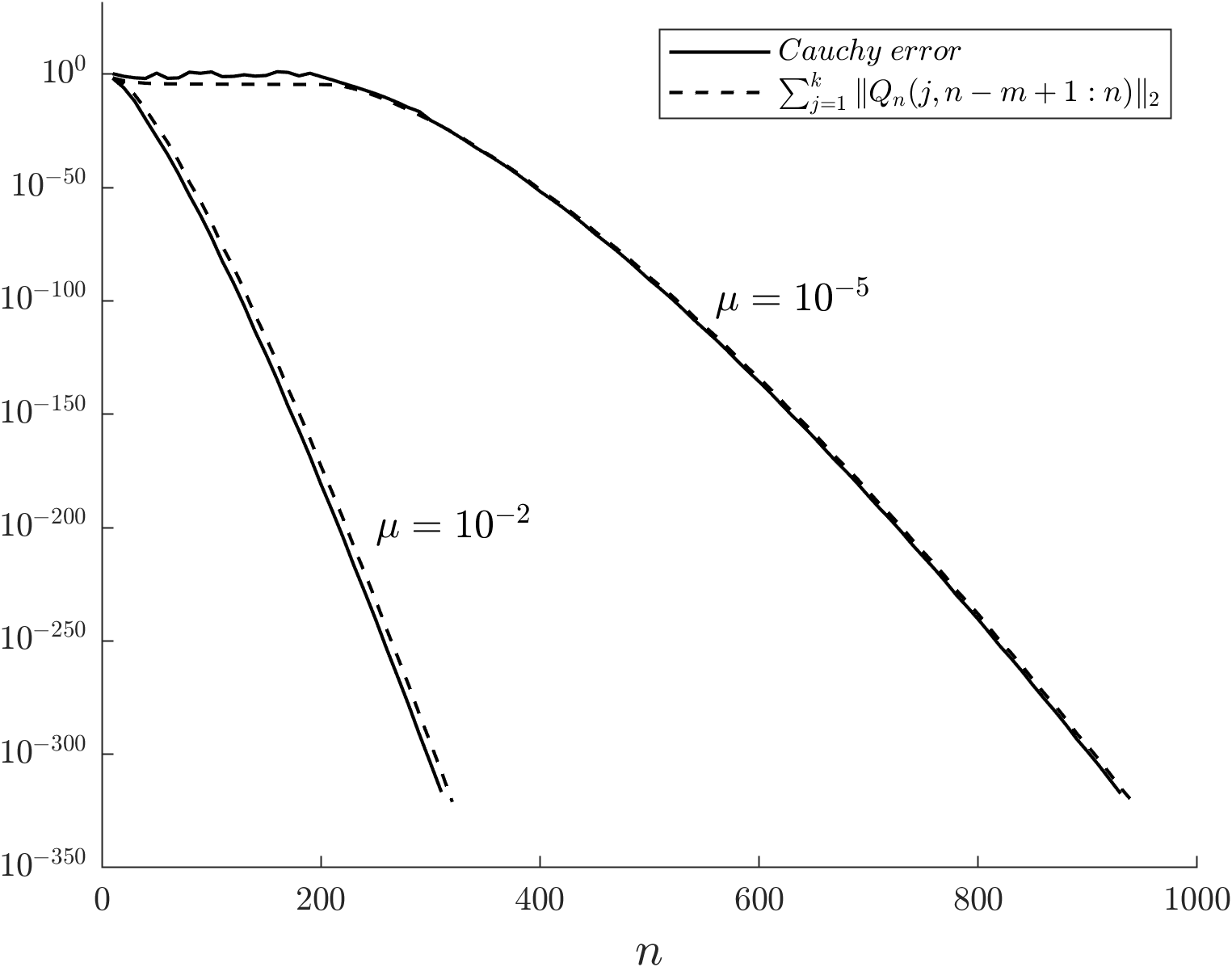}}
  \hfill
  \subfloat[total error]{\label{fig:eps25}\includegraphics[scale=0.46]{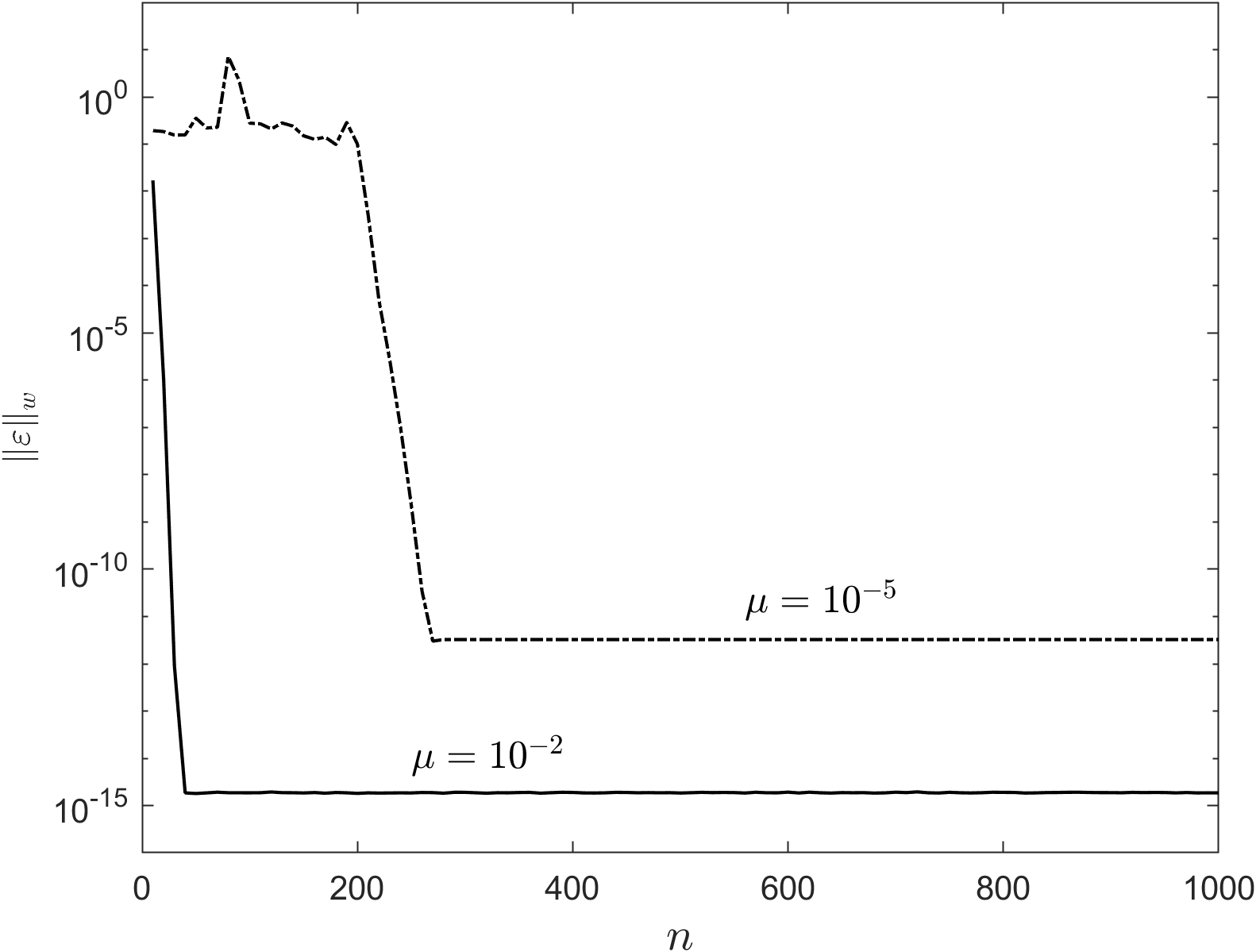}}
  \caption{The explanations of \cref{lem:q} and \cref{lem:u2e0}.}\label{fig:CauchyEps}
\end{figure}

The decay of the Cauchy error is often seen to be exponential or super-exponential, such as the one shown in \cref{fig:cauchy} for the computed solution to \cref{airy} for $\mu = 10^{-2}$ and $\mu = 10^{-5}$. This is mainly due to the (super-)exponential decay of $\|Q(j,n-m+1:n)\|_2$ as $n$ grows (see \cref{fig:CauchyQ}), which, in turn, results in that of the right-hand side $\hs_n$. Thus, the solution $\hu_n$ and the Cauchy error must also exhibit (super-)exponential decay, provided that the entries in $\hR_n$ do not grow (super-)exponentially. For the Airy equation, the decay of $\|Q(j,n-m+1:n)\|_2$ and the Cauchy error match almost perfectly.



A final question we may ask is whether the decaying Cauchy error serves as an indicator of convergence before it reaches the machine epsilon. The answer is yes, owing to Cauchy's limit theorem. However, there is no guarantee that the computed solution converges to the exact solution instead of something else. This is described by Trefethen as \emph{the computed solution is not near the solution, but it is nearly a solution (to a slightly perturbed problem)} \cite[\S 3]{tre2}. When the problem is ill-conditioned, i.e., sensitive to perturbations, the computed solution is very likely to be far off. This is exactly what we see from \cref{fig:eps25}, where the total error $\veps$ of the computed solution to \cref{airy} is plotted for $\mu = 10^{-2}$ and $\mu = 10^{-5}$. The total error plateaus at approximately $3\times 10^{-12}$ for $\mu = 10^{-5}$, suggesting that the decay of the Cauchy error beyond $\mO(10^{-12})$ is worthless despite it signals convergence. In this particular example, the large total error is attributed to $\veps_F$. Hence, caution must be exercised for the use of Cauchy errors in judging convergence and accuracy.

\section{Closing remarks}\label{sec:conclusion}
In this article, we itemize the sources of error of the ultraspherical spectral method, give the condition number of the linear system arising from the ultraspherical spectral method, and explain why it bears a constant value by identifying the effective parts of the system. We also explain the behavior of the Cauchy error and discuss its appropriateness in judging the convergence. One may wonder if the analysis can be applied to other spectral methods with the proviso that the linear system is $m$-Hessenberg and only a finite number of the elements of the right-hand side vector are nonzero. If so, does the analysis end up with the same results? 

\begin{figure}[tbhp]
\centering
\subfloat[total error]{\label{fig:sec4error}\includegraphics[scale=0.46]{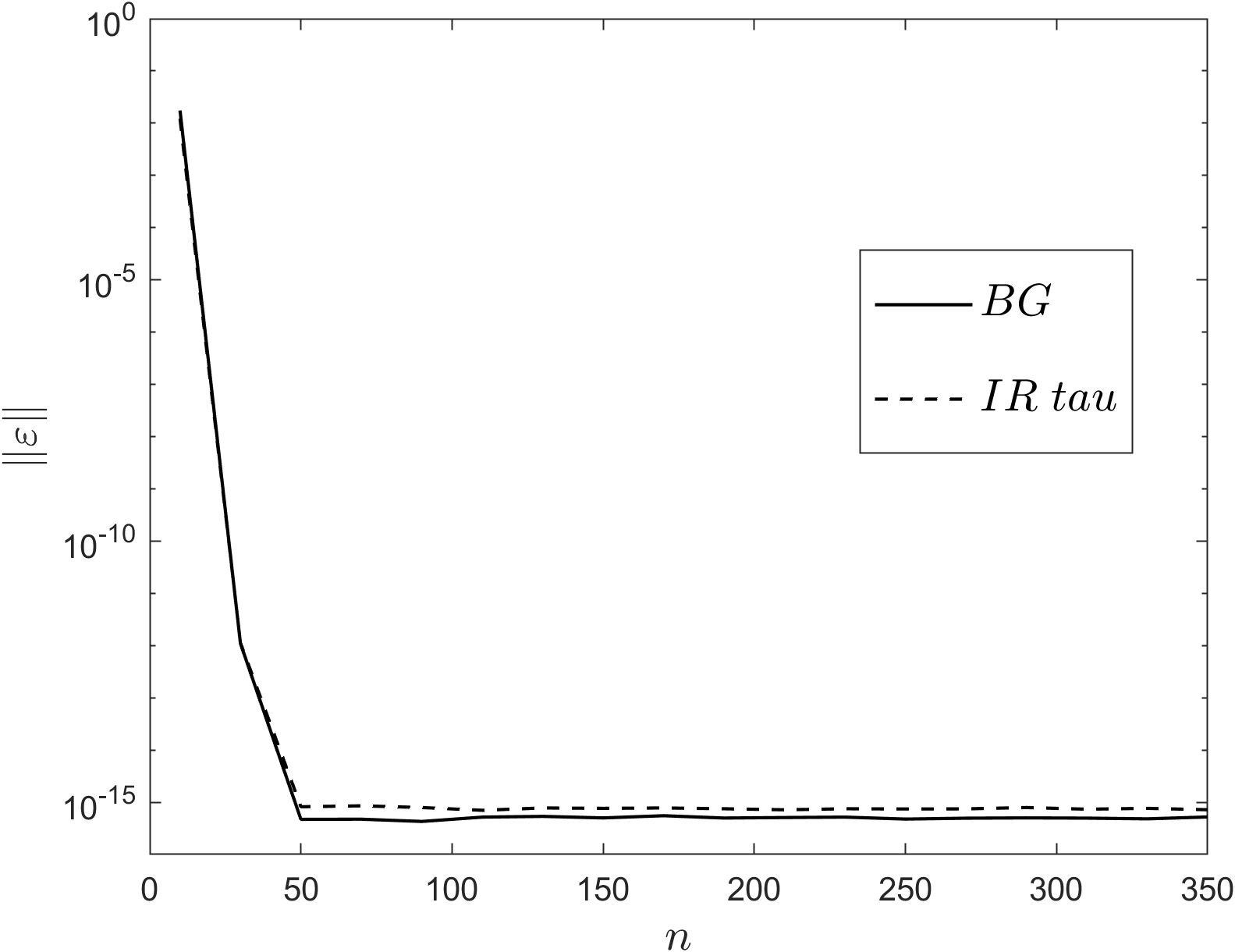}}
\hfill
{\subfloat[Cauchy error]{\label{fig:sec4cauchyerror}\includegraphics[scale=0.46]{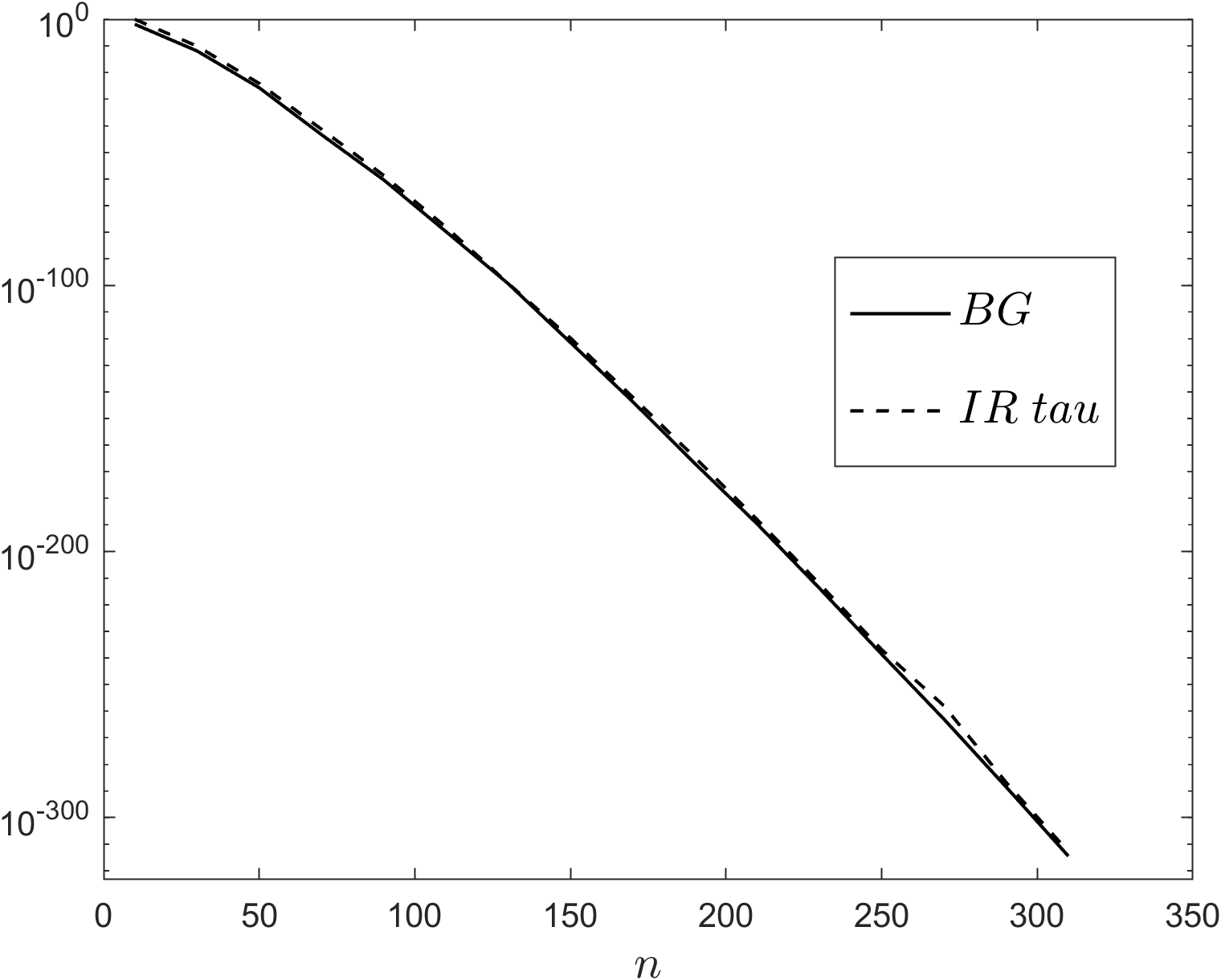}}}
\caption{The total error $\veps$ and the Cauchy error for solving \cref{airy} with $\mu=10^{-2}$ using the banded Galerkin spectral method (BG) and the integral reformulated tau method (IR tau).}\label{fig:sec4}
\end{figure}

We find the answers to these questions affirmative. For example, the banded Galerkin spectral method \cite{mor} and the integral reformulated tau method \cite{du} both yield $m$-Hessenberg matrices and decaying right-hand sides. We solve \cref{airy} for $\mu=10^{-2}$ using these methods and show the results in \cref{fig:sec4}. The total error $\veps$ for both methods are given in \cref{fig:sec4error}, where the plateaus are clear indications of constant condition numbers as in the ultraspherical spectral method. \cref{fig:sec4cauchyerror} shows the Cauchy errors which decay beyond the machine epsilon as expected.

Another natural extension of this work is to partial differential equations. The ultraspherical spectral method for PDEs in multiple dimensions also enjoys constant condition numbers. This can be obtained following an analogous analysis which we omit here. 

\bibliographystyle{siam}
\bibliography{references}

\begin{thebibliography}{1}

\bibitem{aur}
{\sc J.~L. Aurentz and L.~N. Trefethen}, {\em {Chopping a {C}hebyshev series}},
  ACM Transactions on Mathematical Software, 43 (2017), pp.~1--21.

\bibitem{du}
{\sc K.~Du}, {\em On well-conditioned spectral collocation and spectral methods
  by the integral reformulation}, SIAM Journal on Scientific Computing, 38
  (2016), pp.~A3247--A3263.

\bibitem{hig}
{\sc N.~J. Higham}, {\em Accuracy and {S}tability of {N}umerical {A}lorithms},
  vol.~61, SIAM, 1998.

\bibitem{mor}
{\sc M.~Mortensen}, {\em {A generic and strictly banded spectral
  Petrov--Galerkin method for differential equations with polynomial
  coefficients}}, SIAM Journal on Scientific Computing, 45 (2023),
  pp.~A123--A146.

\bibitem{olv}
{\sc S.~Olver and A.~Townsend}, {\em A fast and well-conditioned spectral
  method}, SIAM Review, 55 (2013), pp.~462--489.

\bibitem{tre}
{\sc L.~N. Trefethen}, {\em Approximation {T}heory and {A}pproximation
  {P}ractice, {E}xtended {E}dition}, SIAM, 2019.

\bibitem{tre2}
\leavevmode\vrule height 2pt depth -1.6pt width 23pt, {\em Eight perspectives
  on the exponentially ill-conditioned equation $\varepsilon y''-xy'+y=0$},
  SIAM Review, 62 (2020), pp.~439--462.

\bibitem{wat}
{\sc D.~S. Watkins}, {\em The Matrix Eigenvalue Problem: GR and Krylov Subspace
  Methods}, SIAM, 2007.

\end{thebibliography}

\end{document}